\documentclass[12pt]{amsart}

\usepackage{amsmath,amsthm,amssymb}
\font\teneufm=eufm10
\font\seveneufm=eufm7
\font\fiveeufm=eufm5
\newfam\frakturfam

\textfont\frakturfam=\teneufm
\scriptfont\frakturfam=\seveneufm
\scriptscriptfont\frakturfam=\fiveeufm

\newtheorem{lemma}{Lemma}
\newtheorem{theorem}{Theorem}
\newtheorem{corollary}{Corollary}

\def\bee{\begin{eqnarray}}
\def\bes{\begin{eqnarray*}}
\def\eee{\end{eqnarray}}
\def\ees{\end{eqnarray*}}

\def\a{\alpha}
\def\b{\beta}
\def\g{\gamma}

\def\s{\sigma}
\def\t{\tau}
\def\d{\partial}
\def\Proof{{\sl Proof.}\ }

\title{Automorphisms of free metabelian Lie algebras, I}

\begin{document}
\date{}
\maketitle

\begin{center}

{\bf Ualbai Umirbaev}\footnote{Department of Mathematics,
 Wayne State University,
Detroit, MI 48202, USA; Department of Mathematics, 
Al-Farabi Kazakh National University, Almaty, 050040, Kazakhstan; 
and Institute of Mathematics and Mathematical Modeling, Almaty, 050010, Kazakhstan,
e-mail: {\em umirbaev@wayne.edu}}

\bigskip

\dedicatory{In memory of Professor V.A. Roman'kov (1948--2023)}

\end{center}

\begin{abstract} We show that all Chein automorphisms (or one-row transformations) of lower degree $\geq 4$ of a free metabelian Lie algebra $M_n$ of rank $n\geq 4$ over an arbitrary field $K$ of characteristic $\neq 3$ are tame. We then show that all exponential automorphisms of $M_n$ of lower degree $\geq 5$ are also tame under the same conditions. The same results hold for fields of any characteristic when $n\geq 5$.  These results contradict some long-standing results in the area.

We also prove that a large class of automorphisms of $M_n$ of rank $n\geq 4$ that move only two variables are almost tame, that is, they can be expressed as a product of Chein automorphisms. 

\end{abstract}

\noindent {\bf Mathematics Subject Classification (2020):} 17B40, 17B01, 17B30, 16W20.

\noindent {\bf Key words:} Automorphism, derivation, free metabelian Lie algebra. 

%\maketitle
\tableofcontents

\section{Introduction}

\hspace*{\parindent}

The well known Jung--van der Kulk Theorem \cite{Jung, Kulk} states that every automorphism of the polynomial algebra $K[x,y]$ in two variables $x,y$ over an arbitrary field $K$ is tame. An analogue of this result for free associative algebras in two variables was proven by   Czerniakiewicz and Makar-Limanov \cite{Czer, ML70}. The same result is true for free Poisson algebras in two variables over a field of characteristic zero \cite{MLTU} (see also \cite{MLSh,MLU16,U12}).  Moreover, the automorphism groups of polynomial algebras, free associative algebras, and free Poisson algebras in two variables are isomorphic. The tameness of automorphisms of free right-symmetric algebras in two variables was proven in \cite{KMLU}. 

For any free algebra $A$ in the free variables $x_1,\ldots,x_n$ denote by $\phi=(f_1,\ldots,f_n)$ the endomorphism of $A$ determined by $\phi(x_i)=f_i$ for all $1\leq i\leq n$. 

The automorphism  groups of commutative and associative algebras generated by three elements are much more complicated. 
The well-known Nagata automorphism \cite{Nagata} 
\bes
(x+2y(zx-y^2)+z(zx-y^2)^2, y+z(zx-y^2),z)
\ees
of the polynomial algebra $K[x,y,z]$ over a field $K$ of characteristic $0$ is proven to be wild \cite{SU04-1,US02RAN}. The well-known Anick automorphism (see \cite[p. 398]{Cohn06})
\bes
(x+z(xz-zy),\,y+(xz-zy)z,\,z)
\ees
of the free associative algebra $K\langle x,y,z\rangle$ over a field $K$ of characteristic $0$ is also proven to be wild \cite{UU06AN,UU07JR}. 
The Nagata
automorphism gives an example of a wild
automorphism of free Poisson algebras in three variables. Recently Shestakov and Zhang \cite{SZ24} constructed an analogue of the Anick automorphism for free Poisson algebras in three variables. 

 It is well known \cite{Smith} that the Nagata and Anick automorphisms are stably tame, that is, they become tame after adding one more variable.

A classical theorem of combinatorial group theory asserts that every subgroup of a free group is free. Nielsen established this remarkable result in 1921 for finitely generated subgroups \cite{Nielsen21};  in 1927, it was extended to all subgroups by Schreier \cite{Schreier}. 
 A variety of universal algebras is called Nielsen-Schreier, if any subalgebra of a free algebra of this variety is free. In 1953 Shirshiv and in 1956 Witt proved that  \cite{Shirshov53,Witt} the varieties of Lie algebras and Lie $p$-algebras over a field are Nielsen-Schreier. The varieties of all non-associative algebras \cite{Kurosh},
 commutative and anti-commutative algebras \cite{Shirshov54}, Lie superalgebras  \cite{Mikhalev85,Stern}, and Lie $p$-superalgebras \cite{Mikhalev88} over a field are also Nielsen-Schreier. It was recently shown \cite{DU22} that the varieties of pre-Lie (also known as right-symmetric) algebras and Lie-admissible algebras over a field of characteristic zero are Nielsen-Schreier. 
 Some other examples of Nielsen-Schreier 
 varieties can be found in \cite{Chibrikov,MS14,SU02,U94,U96}.

In 1924, Nielsen  also proved \cite{Nielsen24} that the automorphism group of a finitely generated free group is generated by Nielsen transformations, and represented this group as a finitely defined group. 
 In 1964 P. Cohn  proved \cite{Cohn64} that all automorphisms of finitely generated free
Lie algebras over a field are tame. Later this result was extended to free algebras of Nielsen-Schreier varieties \cite{Lewin}. 
 Defining relations between the elementary automorphisms of free algebras of Nielsen-Schreier varieties are given in \cite{UU07JA}.

One of the classical and well-studied varieties of groups is the variety of metabelian groups. In 1965 Bachmuth \cite{Bachmuth} proved that all automorphisms of the free metabelian group of rank $2$ are tame. The first example of a non-tame automorphism of the free metabelian group of rank $3$ was constructed by Chein \cite{Chein}. 
Bachmuth-Mochizuki proved \cite{BM67} that the group of automorphisms of the free metabelian group of rank $3$ is not even finitely generated. Strengthening these results, Roman'kov proved \cite{Romankov92} that the free metabelian group of rank $3$ contains non-tame primitive elements and  gave a method of constructing of such elements. The well known Bachmuth-Mochizuki-Roman'kov Theorem \cite{BM85,Romankov85} states that every automorphism of the free metabelian group of rank $n\geq 4$ is tame. More information about automorphisms of free metabelian groups can be found in \cite{Romankov20}.

The variety of metabelian Lie algebras is also a classical and well-studied variety of Lie algebras.
Throughout this paper, $M_n$ denotes the free metabelian Lie algebra of rank $n$ over a field $K$ generated by $x_1,\ldots,x_n$. 
It is well known that any nontrivial exponential automorphism of $M_2$ is wild (see, for example \cite{Artamonov,Shmelkin73}). 
In 1992, Drensky \cite{Drensky92} proved that the exponential automorphism $\mathrm{exp}(\mathrm{ad}[x_1,x_2])$ of $M_3$ is wild. More examples of wild automorphisms of $M_3$ are given in \cite{KR,Nauryzbaev09,Romankov08}. 

In 1992 Bahturin and Nabiyev \cite{BN} claimed that every nontrivial exponential automorphism of the free metabelian Lie algebra $M_n$ of any rank $n\geq 2$ over a field of characteristic zero is wild.  More examples of wild automorphisms of $M_n$ of rank $n\geq 4$ were given in 2008 by \"Ozcurt and Ekici \cite{OE}. For a long time, the area specialists believed these results were true. Unfortunately, both articles have fatal errors \cite{U24}.

We show that all Chein automorphisms (or one-row transformations) of lower degree $\geq 4$ of a free metabelian Lie algebra $M_n$ of rank $n\geq 4$ over an arbitrary field $K$ of characteristic $\neq 3$ are tame. We then show that all exponential automorphisms of $M_n$ of lower degree $\geq 5$ are also tame under the same conditions. The same results hold for fields of any characteristic when $n\geq 5$.  These results contradict those of \cite{BN,OE}.

We also prove that a large important class of automorphisms of $M_n$ of rank $n\geq 4$ that move only two variables are almost tame, that is, they can be expressed as a product of Chein automorphisms. 

More precise formulations of the results are given in the following sections.

The paper is organized as follows. In Section 2 we define some important automorphisms of $M_n$ and give precise formulations some of our results.  In Section 3, using a representation of $M_n$ via Magnus embedding and Fox derivatives, we describe Chein autoimorphisms of $M_n$.  In Section 4 we prove that a huge part of Chein and exponential automorphisms are tame. In Section 5 we prove that a huge class of automorphisms that move only two variables are almost tame.

\section{Automorphisms of $M_n$}

\hspace*{\parindent}

Let $M_n$ be the free metabelian Lie algebra over a field $K$ with free generators $x_1,\ldots,x_n$. 
The set of all right normed Lie monomials of the form 
\bee\label{f1}
[\ldots[[x_i,x_{i_1}],x_{i_2}],\ldots,x_{i_k}], \ i>i_1\leq i_2\leq\ldots \leq i_k, 
\eee
is a linear basis of $M_n'=[M_n,M_n]$ (see, for example \cite[p. 139]{Bahturin}). 

Let $\mathrm{Aut}(M_n)$ be the group of all automorphisms of $M_n$. 
   An automorphism of $M_n$ of the form
\bes
(x_{1}, x_{2}, \ldots, x_{i-1}, \alpha x_{i} +f,x_{i+1}, \ldots x_{n} ),
\ees
 where $0 \neq \alpha \in K$ and $f$ does not contain $x_i$, is called {\em elementary}.  The subgroup $\mathrm{TAut}(M_n)$ of $\mathrm{Aut}(M_n)$ generated by all 
  elementary automorphisms is called the {\em tame automorphism group}, 
 and the elements of this subgroup are called {\em tame automorphisms} 
 of $M_n$. Nontame automorphisms of $M_n$ are called {\em wild}. 

Defining relations among elementary automorphisms have been described for polynomial algebras in three variables  \cite{UU06AN,UU06JR}, for finitely generated free Lie algebras \cite{UU07JA},  and for free metabelian Lie algebras in three variables \cite{Nauryzbaev10}.

As usual, the group of all linear automorphisms of $M_n$ will be identified with the group of all invertible matrices $GL_n(K)$. 
For every pair of integers $s,t$ with $1\leq s\neq t\leq n$, let $(st)$ denote the automorphism that simply permutes the variables $x_s$ 
  and $x_t$ and leaves all other variables unchanged. The group generated by all $(st)$ is isomorphic to the symmetric group $S_n$.

In 2015, Nauryzbaev \cite{Nauryzbaev15} proved that the group of tame automorphisms of $M_n$ is generated by all linear automorphisms $GL_n(K)$ and the quadratic automorphism 
\bee\label{f3}
(x_1+[x_2,x_3],x_2,\ldots,x_n)
\eee
if either $n\geq 4$ and $K$ is a field of characteristic $\neq 3$, or $n\geq 5$ and $K$ is a field of arbitrary characteristic.

 If $\phi=(f_1,\ldots,f_n), \psi=(g_1,\ldots,g_n)\in \mathrm{Aut}(M_n)$, then the composition $\phi\psi$ is defined as 
\bes
\phi\psi=(g_1(f_1,\ldots,f_n),\ldots,g_n(f_1,\ldots,f_n)). 
\ees
We also adopt the standard notation $[\phi,\psi]=\phi\psi\phi^{-1}\psi^{-1}$ for the commutator, and $\phi^{\psi}=\psi\phi\psi^{-1}$ for conjugation.

Let $0\neq z\in [M_n,M_n]$. Then $\mathrm{ad}\,z : M_n\to M_n (m\mapsto [z,m])$ is a nilpotent derivation and $(\mathrm{ad}\,z)^2=0$. The {\em exponential} automorphism 
$\exp(\mathrm{ad}\,z)$ of $M_n$ can be written as  
\bes
\exp(\mathrm{ad}\,z)=\mathrm{id}+\mathrm{ad}\,z =(x_1+[z,x_1], \ldots, x_n+[z,x_n]). 
\ees
 It is known that the exponential automorphism $\exp(\mathrm{ad}\,[x_1,x_2])$ of $M_3$ is wild \cite{Drensky92}.

Recall that in 1993 Bryant and Drensky proved \cite{BD93-2} that every automorphism of $M_n$ of rank $n\geq 4$ over a field $K$ of characteristic zero is {\em approximately tame} ( see the definitions in \cite{SU25}); that is, it can be approximated by tame automorphisms with respect to power series topology \cite{BD93-2,SU25}. In contrast, all wild automorphisms of $M_3$ are {\em absolutely wild} \cite{SU25}, meaning that they cannot be approximated by tame automorphisms.

Free metabelian Lie algebras admit automorphisms of the form  
\bes
(x_1,\ldots,x_{i-1},f_i,x_{i+1},\ldots,x_n), \ 1\leq i\leq n, 
\ees
which include all elementary automorphisms but not necessarily elementary ones. They are called {\em Chein automorphisms} or {\em one-row transformations} \cite{GGR,Romankov95}. These automorphisms move only one variable and fix the others. The subgroup $\mathrm{ATAut}(M_n)$ of $\mathrm{Aut}(M_n)$ generated by all 
  Chein automorphisms is called the {\em almost tame automorphism group}, 
 and the elements of this subgroup is called {\em almost tame automorphisms} 
 of $M_n$ \cite{Nauryzbaev09,Nauryzbaev11}. Thus 
\bes
\mathrm{TAut}(M_n)\subseteq \mathrm{ATAut}(M_n)\subseteq \mathrm{Aut}(M_n). 
\ees

Obviously, the group of almost tame automorphisms $\mathrm{ATAut}(M_n)$ is generated by all linear automorphisms $GL_n(K)$ and all Chein automorphisms of the form  
\bee\label{f4}
(x_1+f, x_2,\ldots,x_n), \ f\in [M_n,M_n]. 
\eee
Chein automorphisms of this form are described in Section 3.

 The lowest-degree example of a non-elementary Chein automorphism is
\bee\label{f5}
(x_1+[[x_2,x_3],x_1],x_2,\ldots,x_n). 
\eee
This automorphism is wild if $n=3$ \cite{Nauryzbaev09,Romankov08}, and is the first candidate for wildness for all $n\geq 4$, but it is still unknown whether it is tame or not.  It was recently shown \cite{NSU25} that an analogue of this automorphism for free metabelian anticommutative algebras is absolutely wild for all $n\geq 3$. 

In 2008, Roman'kov \cite{Romankov08} proved that the group of automorphisms of $M_3$ cannot be finitely generated modulo the subgroup generated by all exponential and Chein automorphisms. Kabanov and Roman'kov \cite{KR} constructed strictly nontame primitive elements of $M_3$. Defining relations  for Chein automorphisms of $M_3$ are given in \cite{Nauryzbaev11}.

We introduce the notion of the lower degree of an automorphism. Consider the standard degree grading of 
\bes
M_n=L_1\oplus L_2\oplus \ldots\oplus L_k\oplus\ldots, 
\ees
where $L_k$ is the linear span of all Lie monomials of degree $k$ in $x_1,\ldots,x_n$. Every nonzero $m\in M_n$ can be uniquely written in the form 
\bes
m=l_s+l_{s+1}+\ldots+l_t, 
\ees
where $l_j\in L_j$ for all $s\leq j\leq t$, $l_s\neq 0$, and $l_t\neq 0$.  We call $\mathrm{ldeg}(m)=s$ the {\em lower degree}  of $m$, while $\deg(m)=t$ denotes the degree of $m$. If $\mathrm{ldeg}(m)=\deg(m)=t$, then $m$ is {\em homogeneous} of degree $t$. Set also $\mathrm{ldeg}(0)=\infty$ and $\deg(0)=-\infty$. 

Let 
\bes
\phi=(l_1+m_1,\ldots,l_n+m_n)
\ees
be an automorphism of $M_n$ such that $l_i\in L_1$ and $m_i\in [M_n,M_n]$ for all $i$. Set 
\bes
\mathrm{ldeg}(\phi)=\mathrm{min}\{\mathrm{ldeg}(m_i) | 1\leq i\leq n\}
\ees
and 
\bes
\mathrm{deg}(\phi)=\mathrm{max}\{\mathrm{deg}(m_i)| 1\leq i\leq n\}. 
\ees
If $\mathrm{ldeg}(\phi)=\mathrm{deg}(\phi)=t$, then we say that $\phi$ is {\em homogeneous} of degree $t$. In particular, (\ref{f3}) is quadratic and (\ref{f5}) is cubic. 
Clearly, every nonlinear automorphism has a finite lower degree and degree. 

We show that all Chein automorphisms of lower degree $\geq 4$ are tame. This means that all Chein automorphisms "except" (\ref{f5}) are tame. We also show that all exponential automorphisms of lower degree $\geq 5$ are tame. These results contradict the results of \cite{BN,OE}.

\section{Chein automorphisms}

\hspace*{\parindent} 

In this section, using an analogue of the well known Magnus embedding \cite{KSh,Magnus}, we give another representation of the free metabelian algebra $M_n$, and give direct definition of the Fox derivatives for $M_n$. This construction for Lie algebras is a particular case of Shmel'kin's  wreath products \cite{Shmelkin73}. 

Let $Y_n=M_n/[M_n,M_n]$ be the abelian Lie algebra with a linear basis $y_1,y_2,\ldots,y_n$, where $y_i=x_i+[M_n,M_n]$. 
The universal enveloping algebra $U(Y_n)$ of $Y_n$ is the polynomial algebra $U=U(Y_n)=K[y_1,\ldots,y_n]$ in the variables $y_1,y_2,\ldots,y_n$. Let 
\bes
U=U_0\oplus U_1\oplus \ldots \oplus U_k\oplus \ldots
\ees
be the degree grading of $U$, i.e., $U_k$ is the linear span of all homogeneous monomials in $y_1,y_2,\ldots,y_n$ of degree $k$.
If 
\bes
a=a_s+a_{s+1}+\ldots+a_t\in U, a_j\in U_j, a_s\neq 0, a_t\neq 0, 
\ees
then $\deg(a)=t$ is the degree and $\mathrm{ldeg}(a)=s$ is the lower degree of $a$ in the variables $y_1,y_2,\ldots,y_n$.
Set $\deg(0)=-\infty$ and $\mathrm{ldeg}(0)=\infty$ as usual.

Let $T_n$ be the free right $U$-module with basis $t_1,\ldots,t_n$. Turn the direct sum 
\bes
M=Y_n\oplus T_n 
\ees
  into a Lie algebra by 
\bes
[a+t,b+s]= tb-sa, 
\ees
where $a,b\in Y_n$ and $t, s\in T_n$. Then the subalgebra of $M$ generated by 
\bes
x_i=y_i+t_i, \ \ 1\leq i\leq n, 
\ees
is the free metabelian Lie algebra  with free generators $x_1,\ldots,x_n$ \cite{KSh,Magnus,Shmelkin73}. We identify this algebra with $M_n$. 

Notice that $M_n'=[M_n,M_n]\subseteq T_n$ and for any $m\in M_n'$ and $y_i$ we get 
\bes
m y_i=[m,x_i]. 
\ees

Every element $m\in M_n'$ can be written in the form 
\bes
m=\sum_{1\leq i<j\leq n} [x_i,x_j]a_{ij}, \ \ a_{ij}\in U. 
\ees

Let $f\in M_n$ be an arbitrary element $M_n$. Then $f$ is uniquely represented as 
\bes
f=y+t, y\in Y_n, t\in T_n. 
\ees
Moreover, $t$ is uniquely represented as 
\bes
t=t_1d_1+\ldots+t_nd_n, \ \ \ \ d_1,\ldots,d_n\in K[y_1,\ldots,y_n]=U. 
\ees
Set  $d_i=\frac{\partial f}{\partial x_i}$ for all $1\leq i\leq n$ and 
\bes
\partial(f)=(d_1,\ldots,d_n)^t=\big(\frac{\partial f}{\partial x_1},\ldots,\frac{\partial f}{\partial x_n}\big)^t, 
\ees
where $t$ means the transpose. We call $\frac{\partial f}{\partial x_1},\ldots,\frac{\partial f}{\partial x_n}$ the Fox derivatives of $f$  (see \cite{Shpilrain,U93}). 

Let $e_1,e_2,\ldots,e_n$ be the columns of the identity matrix $I=I_n$ of order $n$. 

The statement of the next lemma is well known \cite{U93,95SMJ}. 
\begin{lemma}\label{l1}  Let $a=(a_1,a_2,\ldots,a_n)^t\in (U^n)^t$ be an arbitrary column of dimension $n$ over $U$. 
Then the following conditions are equivalent. 

$(a)$ $a=\d(f)$ for some $f\in [M_n,M_n]$; 

$(b)$  $Ya=y_1a_1+\ldots +y_na_n=0$ where $Y=(y_1,\ldots,y_n)$; 

$(c)$ $a=\sum_{i<j} (e_iy_j-e_jy_i)a_{ij}, \ a_{ij}\in U$. 
\end{lemma}
\Proof 
If $i<j$ and $a\in U$, then 
\bes
\d([x_i,x_j]a)=e_iy_ja-e_jy_ia=(e_iy_j-e_jy_i)a. 
\ees
Consequently, $(a)$ implies $(c)$. Obviously, $(c)$ implies $(b)$ since $Y(e_iy_j-e_jy_i)=0$ for all $i<j$. 

Suppose that $(b)$ holds. Let $a_i=-y_nb_i+a_i^0$, where $a_i^0\in K[y_1,\ldots,y_{n-1}]$ for all $1\leq i\leq n-1$. 
Then  $(b)$ implies that $a_n=y_1b_1+\ldots+y_{n-1}b_{n-1}$. Set $g=[x_n,x_1]b_1+\ldots+[x_n,x_{n-1}]b_{n-1}$. We get  
$\frac{\d g}{\d x_n}=a_n$ and the last component of $a'=a-\d(g)$ is zero. Moreover, if $a'=(a_1',\ldots,a_{n-1}',0)$ then
$y_1a_1'+\ldots+y_{n-1}a_{n-1}'=0$. Leading an induction on the number of nonzero components of $a$, we may assume that there exists $h\in [M_n,M_n]$ such that $a'=\d(h)$. Then $a=\d(g+h)$. 
$\Box$

Every endomorphism $\phi$ of $M_n$ induces an endomorphism of $Y_n=M_n/[M_n,M_n]$, which can be uniquely extended to an endomorphism of $U=U(Y_n)$. This endomorphism will be denoted again by $\phi$ for simplicity of notation. If $f\in U$, then the image $\phi(f)$ will be denoted by $f^{\phi}$ for convenience.  

We define the Jacobian matrix of $\phi=(f_1,\ldots,f_n)$  by 
\bes
J(\phi)=\left[\frac{\partial f_j}{\partial x_i}\right]_{1\leq i,j\leq n}=\left[\partial(f_1) \  \ldots \ \partial(f_n)\right].
\ees
If $\phi=(f_1,\ldots,f_n)$ and $\psi=(g_1,\ldots,g_n)$ then 
\bes
\phi\psi(x_i)=g_i(f_1,\ldots,f_n)
\ees
for all $1\leq i\leq n$. The corresponding Chain Rule gives \cite{Shpilrain,U93,95SMJ} that 
\bes
J(\phi\psi)=J(\phi) J(\psi)^{\phi}. 
\ees
Consequently, if $\phi$ is an automorphism, then $J(\phi)$ is invertible. The converse is also true  \cite{Shpilrain,U93}. 
\begin{lemma}\label{l2} 
An endomorphism $\phi$ of $M_n$ is an automorphism if and only if the Jacobian matrix $J(\phi)$ is invertible.  
\end{lemma}
\Proof Suppose that $J(\phi)$ is invertible. Without loss of generality, we may assume that 
\bes
\phi=(x_1+f_1,\ldots,x_n+f_n), \ \ f_1,\ldots,f_n\in [M_n,M_n]. 
\ees
Then 
\bes
J(\phi)=I+A, \ \ A=[\d(f_1),\ldots,\d(f_n)]. 
\ees
We get  
\bes
J(\phi)^{-1}=I+B, B=\sum_{i\geq 1} (-A)^i. 
\ees
By Lemma \ref{l1}, we have $YA=0$. Consequently, $YB=0$. Again, by Lemma \ref{l1}, there exist $g_1,\ldots,g_n\in [M_n,M_n]$ such that 
\bes
B=[\d(g_1),\ldots,\d(g_n)]. 
\ees
 Let $\psi=(x_1+g_1,\ldots,x_n+g_n)$. Then 
\bes
J(\psi)=I+B=J(\phi)^{-1}
\ees
 and 
\bes
J(\phi\psi)=J(\phi)J(\psi)^{\phi}= J(\phi)J(\psi)=I. 
\ees
Consequently, $\phi\psi=\mathrm{Id}$ and $\phi$ is an automorphism. $\Box$

We now give a description of Chein automorphisms of the form (\ref{f4}). 

\begin{lemma}\label{l4} An endomorphism of $M_n$ of the form (\ref{f4}) 
 is an automorphism if and only if $f$ belongs to the ideal generated by all $[x_s,x_t]$, where $s,t\neq 1$.   
\end{lemma}
\Proof Let $\phi=(x_1+f,x_2,\ldots,x_n)$ where $f\in [M_n,M_n]$. The Jacobian matrix of $\phi$ is 
\bes
J(\phi)=
I+[\d(f) \ 0\ldots 0]. 
\ees 
This matrix is invertible if and only if $1+\frac{\d f}{\d x_1}$ is invertible in $U$ or, equivalently, 
$\frac{\d f}{\d x_1}=0$. We want to show that $\frac{\d f}{\d x_1}=0$ if and only if $f$ belongs to the ideal generated by all $[x_s,x_t]$, where $s,t\neq 1$. 

The basis (\ref{f1}) of $[M_n,M_n]$ is convenient to work with $\frac{\d f}{\d x_n}=0$. Any $f\in [M_n,M_n]$ can be written in the form 
\bes
f=f'+[x_n,x_1]a_1+[x_n,x_2]a_2+\ldots+[x_n,x_{n-1}]a_{n-1}, 
\ees
where $f'$ belongs to the ideal generated by all $[x_s,x_t]$ with $s,t\neq n$ and $a_i\in K[y_i,\ldots,y_n]$ for all $i$. 
Since $\frac{\d f'}{\d x_n}=0$ it follows that 
\bes
\frac{\d f}{\d x_n}=y_1a_1+y_2a_2+\ldots+y_{n-1}a_{n-1}. 
\ees
If 
\bes
y_1a_1+y_2a_2+\ldots+y_{n-1}a_{n-1}=0
\ees
then the substitution $y_1=\ldots=y_{n-2}=0$ gives that $a_{n-1}=0$. Repeated substitution $y_1=\ldots=y_{n-3}=0$ gives $a_{n-2}=0$. 
Continuing the same discussions, we get $a_{n-1}=a_{n-2}=\ldots=a_1=0$. Therefore $\frac{\d f}{\d x_n}=0$ if and only if $f=f'$. 
$\Box$

\section{Tame automorphisms}

\hspace*{\parindent} 

In this section we always assume that either $n\geq 4$ and $K$ is a field of characteristic $\neq 3$ or $n\geq 5$.

\begin{theorem}\label{t1} Every Chein automorphism of $M_n$ of lower degree $\geq 4$ is tame. 
\end{theorem}
\Proof It suffices to prove that every Chein automorphism of the form 
\bes
C(a)= (x_1+[x_2,x_3]a, x_2,\ldots,x_n), \ a\in U, 
\ees
with $\mathrm{ldeg}(a)\geq 2$ is tame. Since
\bes
C(a+b)=C(a)C(b), \ \ a,b\in U, 
\ees
we may assume that $a=\g y_1^{i_1}y_2^{i_2}\ldots y_n^{i_n}$, $0\neq \g\in K$. If $i_1=0$ then $C(a)$ is elementary. Assume that $i_1\geq 1$. 

{\em Case 1.} Suppose that $i_j\geq 1$ for some $2\leq j\leq n-1$. Consider the elementary automorphisms  
\bes
\phi=(x_1,\ldots,x_{n-1}, x_n+[x_2,x_3]y_1^{i_1}\ldots y_j^{i_j-1}),\\
\psi=(x_1-\g[x_j,x_n]y_{j+1}^{i_{j+1}}\ldots y_n^{i_n}, x_2,\ldots,x_n). 
\ees
Direct calculations give that  $[\phi,\psi]=C(a)$ is tame.

{\em Case 2.} Let $i_2=i_3=\ldots=i_{n-1}=0$ and $i_n\neq 0$. Then $a=\g y_1^{i_1}y_n^{i_n}$. We can make $\g=1$ by conjugating 
 $C(a)$ with $(x_1,\g^{-1} x_2,x_3\ldots,x_n)$. Then assume that $a=y_1^{i_1}y_n^{i_n}$. We have 
\bee\label{f8}
[x_2,x_3]a=[x_2,x_n]y_1^{i_1}y_3y_n^{i_n-1} + [x_n,x_3]y_1^{i_1}y_2y_n^{i_n-1} 
\eee
by the Jacobi identity. 
If $i_n\geq 2$ then 
\bes
\phi=(x_1+[x_2,x_n]y_1^{i_1}y_3y_n^{i_n-1},x_2,\ldots,x_n)=(x_1+[x_2,x_3]y_1^{i_1}y_3^{i_n-1}y_n,x_2,\ldots,x_n)^{(3n)}
\ees
and
\bes
\psi=(x_1+[x_n,x_3]y_1^{i_1}y_2y_n^{i_n-1}, x_2,\ldots,x_n)=(x_1+[x_2,x_3]y_1^{i_1}y_2^{i_n-1}y_n, x_2,\ldots,x_n)^{(2n)}
\ees
are tame by Case 1. Therefore $C(a)=\phi\psi$ is tame. 

If $i_n=1$, then set $\a=(x_1,x_2+x_n,x_3,\ldots,x_n)$. We have 
\bes
\tau=(x_1+[x_2,x_3] y_1^{i_1}y_2 ,x_2,\ldots,x_n)
\ees
 is tame by Case 1. Direct calculations give that 
\bes
[\a,\tau]=(x_1+([x_2,x_3]y_n+[x_n,x_3]y_2)y_1^{i_1},x_2,\ldots,x_n)(x_1+[x_n,x_3]y_1^{i_1}y_n,x_2,\ldots,x_n)
\ees
is tame. Since the second factor 
\bes
(x_1+[x_n,x_3]y_1^{i_1}y_n,x_2,\ldots,x_n)=(x_1+[x_2,x_3]y_1^{i_1}y_2,x_2,\ldots,x_n)^{(2n)}=\tau^{(2n)}
\ees
of this product is also tame it follows that its first factor 
\bes
\s=(x_1+([x_2,x_3]y_n+[x_n,x_3]y_2)y_1^{i_1},x_2,\ldots,x_n)
\ees
is tame. 
By (\ref{f8}), we get 
\bes
3[x_2,x_3]a=([x_2,x_3]y_n+[x_n,x_3]y_2)y_1^{i_1}+([x_2,x_3]y_n+[x_2,x_n]y_3)y_1^{i_1}. 
\ees
Consequently, 
\bes
C(3a)=(x_1+([x_2,x_3]y_n+[x_n,x_3]y_2)y_1^{i_1},x_2,\ldots,x_n)\\
\times (x_1+([x_2,x_3]y_n+[x_2,x_n]y_3)y_1^{i_1},x_2,\ldots,x_n)=\s (\s^{-1})^{(23)}
\ees
is tame. If the characteristic of $K$ is not $3$, then its conjugation by $(x_1,3x_2,x_3,\ldots,x_n)$ gives $C(a)$. 

If $n\geq 5$, then conjugating $C(a)$ by $(4n)$, we can make $i_4>0$. Case 1 gives that $C(a)$ is tame again.

{\em Case 3.} Let $i_2=i_3=\ldots=i_{n-1}=i_n=0$. In this case $a=y_1^{i_1}$ and $s=i_1\geq 2$ since $\mathrm{ldeg}(a)\geq 2$.  Case 2 gives that 
\bes
\phi=(x_1+[x_2,x_3]y_1^{s-1}y_n,x_2,\ldots,x_n)
\ees
is tame. Set $\b=(x_1,,\ldots,x_{n-1},x_n+x_1)$. Direct calculations give that  
\bes
[\b,\phi]=C(a) (x_1,\ldots,x_{n-1},x_n-[x_2,x_3]y_1^{s-1}(y_1+y_n)\\
=C(a) (x_1-[x_2,x_3]y_n^s,\ldots,x_n)^{(1n)}(x_1-[x_2,x_3]y_1y_n^{s-1},x_2,\ldots,x_n)^{(1n)}. 
\ees
Since the second factor of this product is tame and the third factor is tame by Case 2 it follows that $C(a)$ is tame.  
$\Box$

This theorem contradicts the results of \cite{OE}.

\begin{corollary}\label{c1} 
The group of all almost tame automorphisms $\mathrm{ATAut}(M_n)$ is generated, modulo all linear automorphisms, by the quadratic automorphism (\ref{f3}) and the  cubic automorphism (\ref{f5}).
\end{corollary}
\Proof The group of all tame automorphisms $\mathrm{TAut}(M_n)$ is generated,  modulo all linear automorphisms, by the quadratic automorphism (\ref{f3}) \cite{Nauryzbaev15}. Obviuosly, the group of all almost tame automorphisms $\mathrm{ATAut}(M_n)$ is generated, modulo all linear automorphisms, by Chein automorphisms $C(a)$. We have 
\bes
C(a)=C(a_0)C(a_1)\ldots C(a_k), 
\ees
where $a_i\in U_i$ for all $0\leq i\leq k$. We have $C(a_0)$ is elementary since $a_0\in K$ and $C(a_2),\ldots,C(a_k)$ are tame by Theorem \ref{t1}. 

Let 
\bes 
a_1=\a_1 y_1+\a_2y_2+\ldots+\a_ny_n.
\ees
 Then 
\bes
C(a_1)=C(\a_1y_1)C(\a_2y_2)\ldots C(\a_ny_n).
\ees
 All $C(\a_2y_2),\ldots, C(\a_ny_n)$ are again tame. 

Consider $C(\a_1y_1)$. Conjugating this by $(x_1,\a_1^{-1}x_2,x_3,\ldots,x_n)$ we can make $\a_1=1$ if $\a_1\neq 0$. 

Consequently, the group of all almost tame automorphisms $\mathrm{ATAut}(M_n)$ modulo $\mathrm{TAut}(M_n)$  is generated by $C(y_1)$. Note that $C(y_1)$ is the cubic automorphism (\ref{f5}).  Therefore the group of all almost tame automorphisms $\mathrm{ATAut}(M_n)$ is generated, modulo all linear automorphisms, by the quadratic automorphism (\ref{f3}) and the cubic automorphism (\ref{f5}). $\Box$

\begin{lemma}\label{l5} Every automorphism of $M_n$ of the form 
\bes
D(a)=(x_1+[x_1,x_2] y_1a, x_2+[x_1,x_2] y_2a,x_3,\ldots,x_n)
\ees
is tame if $\mathrm{ldeg}(a)\geq 2$. 
\end{lemma}
\Proof By Theorem \ref{t1}, if $\mathrm{ldeg}(b)\geq 1$, then 
\bes
\phi=(x_1-[x_3,x_4]y_1b,x_2,\ldots,x_n)(x_1,x_2-[x_3,x_4]y_2b,x_3,\ldots,x_n)\\
=(x_1-[x_3,x_4]y_1b,x_2-[x_3,x_4]y_2b,x_3,\ldots,x_n)
\ees
is tame. Let 
\bes
\psi=(x_1,x_2,x_3+[x_1,x_2],x_4,\ldots,x_n). 
\ees
Then 
\bes
[\phi,\psi]=(x_1+[x_1,x_2]y_1y_4b,x_2+[x_1,x_2]y_2y_4b,x_3,\ldots,x_n)=D(y_4b)
\ees
is tame. 

Consequently, the statement of the lemma is true if $y_4$ divides $a$. Obviously, the same is true if $y_i$ divides $a$ for all $3\leq i\leq n$. Since 
\bes
D(a_1+a_2)=D(a_1)+D(a_2), \ \ a_1,a_2\in U, 
\ees
 we may assume that $a\in K[y_1,y_2]$. Let $\a=(x_1,x_2,x_3,x_4+x_1,x_5,\ldots,x_n)$. 
 If $b\in K[y_1,y_2]$ then
\bes
D(y_4b)^{\a}=(x_1+[x_1,x_2]y_1(y_4+ y_1)b, x_2+[x_1,x_2]y_2(y_4+ y_1)b, \\
x_3, x_4-[x_1,x_2]y_2(y_4+ y_1),x_5,\ldots,x_n)\\
=D(y_4b)D(y_1b)(x_1,x_2,x_3, x_4-[x_1,x_2]y_2(y_4+ y_1),x_5,\ldots,x_n). 
\ees
Since $D(y_4b)$ is tame and 
\bes(x_1,x_2,x_3, x_4-[x_1,x_2]y_2(y_4+ y_1),x_5,\ldots,x_n)
\ees
 is tame by Theorem \ref{t1}, this implies that $D(y_1b)$ is also tame. Therefore the statement of the lemma is true if $a\in K[y_1,y_2]$ and $y_1$ divides $a$. 

Changing $\a$ to $\b=(x_1,x_2,x_3,x_4+x_2,x_5,\ldots,x_n)$ and repeating the same discussions as above, we can show that $D(y_2b)$ is also tame. $\Box$

\begin{theorem}\label{t2}  Every exponential automorphism of $M_n$ of lower degree $\geq 5$ is tame. 
\end{theorem}
\Proof Let 
\bes
E(m)=\exp(\mathrm{ad}(m))=(x_1+[m,x_1], x_2+[m,x_2],\ldots,x_n+[m,x_n])
\ees
be an exponential automorphism with $\mathrm{ldeg}(m)\geq 4$. Since 
\bes
E(m_1+m_2)=E(m_1)+E(m_2), 
\ees
 we may assume that $m=[x_1,x_2]a$ with $\mathrm{ldeg}(a)\geq 2$. Then
\bes
E(m)=E([x_1,x_2]a)=D(a)\phi_3\ldots \phi_n, 
\ees
where 
\bes
\phi_i=(x_1,\ldots,x_{i-1},x_i+[x_1,x_2]y_ia,x_{i+1},\ldots,x_n), \  3\leq i\leq n. 
\ees
All $\phi_i$ are tame by Theorem \ref{t1} and $D(a)$ is tame by Lemma \ref{l5}. 
Consequently, $E(m)$ is tame. $\Box$

This theorem contradicts the main result of \cite{BN}.

\section{Almost tame automorphisms}

\hspace*{\parindent} 

In this section we prove that a huge class of automorphisms of $M_n$ that moves only two variables $x_1$ and $x_2$ are almost tame. We always assume that $n\geq 4$ and $K$ is an arbitrary field.  Since all Chein automorphisms are automatically included into $\mathrm{ATAut}(M_n)$ by the definition, we don't need any restrictions to  characteristic of the field from the preceeding section. 

First we give minor extensions of Lemma \ref{l5} and Theorem \ref{t2} for almost tame automorphisms. 
\begin{lemma}\label{l6} Every automorphism of the form 
\bes
D(a)=(x_1+[x_1,x_2] y_1a, x_2+[x_1,x_2] y_2a,x_3,\ldots,x_n)
\ees
is tame if $\mathrm{ldeg}(a)\geq 1$. 
\end{lemma}
\Proof For any $b\in U$ we have 
\bes
\phi=(x_1-[x_3,x_4]y_1b,x_2,\ldots,x_n)(x_1,x_2-[x_3,x_4]y_2b,x_3,\ldots,x_n)\\
=(x_1-[x_3,x_4]y_1b,x_2-[x_3,x_4]y_2b,x_3,\ldots,x_n)
\ees
is almost tame.
Repeating the proof of Lemma \ref{l5} with this data, we get that $D(a)$ is almost tame for all $a$ with $\mathrm{ldeg}(a)\geq 1$.  $\Box$

\begin{corollary}\label{c2} Every exponential automorphism of lower degree $\geq 4$ is almost tame. 
\end{corollary}
\Proof Consider the decomposition  of the exponential automorphism $E([x_1,x_2]a)$ from the proof of Theorem \ref{t2}.   Note that  $\mathrm{ldeg}(E([x_1,x_2]a))\geq 4$ implies $\mathrm{ldeg}(a)\geq 1$. Consequently, $D(a)$ is almost tame by Lemma \ref{l6}. Since all $\phi_i$ are Chein automorphisms it follows that $E([x_1,x_2]a)$ is almost tame. Consequently, every exponential automorphism of lower degree $\geq 4$ belongs to $\mathrm{ATAut}(M_n)$. $\Box$

\begin{lemma}\label{l7} Every automorphism of the form 
\bes
A(h,g)=(x_1+[x_1,x_n]hg+[x_2,x_n]hg^2, x_2-[x_1,x_n]h-[x_2,x_n]hg,x_3,\ldots,x_n),
\ees
where  $g,h\in U$, is almost tame. 
\end{lemma}
\Proof $(a)$ Let  $\a=(x_1-\lambda x_2, x_2,\ldots,x_n)$. Then 
\bes
A(h,g)^{\a}=A(h^{\a},g^{\a}-\lambda). 
\ees
 This means that $A(h,g)$ is almost tame if and only if $A(h^{\a},g^{\a}-\lambda)$ is almost tame. Using this we may assume that $g(0,\ldots,0)=0$. 

$(b)$ Consider the almost tame automorphisms 
\bes
\phi=(x_1-[x_2,x_3]a,x_2,\ldots,x_n) 
\ees
and 
\bes
\psi=(x_1+[x_3,x_n]y_2ah(g-y_3a),x_2,\ldots,x_n)(x_1,x_2-[x_3,x_n]y_2ah,x_3,\ldots,x_n)\\
=(x_1+[x_3,x_n]y_2ah(g-y_3a),x_2-[x_3,x_n]y_2ah,x_3,\ldots,x_n)
\ees
Direct calculations give that 
\bes
A(h,g)^{\phi}=A(h,g-y_3a)\psi. 
\ees
Consequently, $A(h,g)$ is almost tame if and only if $A(h,g-y_3f)$ is almost tame. Of course, the same trick can be done with any variable $y_i$ if $3\leq i\leq n-1$. 
Therefore, replacing $g$ by $g-y_if$ for all $3\leq i\leq n-1$, we may assume that  
\bes
g=G(y_1,y_2)-y_nb, \ \ b\in U. 
\ees

$(c)$ Now we show how to eliminate $y_nb$. First we replace $g$ to 
\bes
g+y_3b=G(y_1,y_2)+y_3b-y_nb
\ees
using $(b)$. If $\b=(x_1,x_2,x_3+x_n,x_4,\ldots,x_n)$, then it is easy to check that 
\bes
A(h,g+y_3b)^{\b}=A(h^{\b},(g+y_3b)^{\b}). 
\ees
We have 
\begin{align*}
(g+y_3b)^{\b}=g^{\b}+(y_3+y_n)b^{\b}\\
=G(y_1,y_2)-y_nb^{\b}+(y_3+y_n)b^{\b}=G(y_1,y_2)+y_3b^{\b}. 
\end{align*}
Using $(b)$ again, we may assume that 
\bes
g=G(y_1,y_2). 
\ees

(d) Let $c$ be a polynomial in $y_1,y_2$. Then $g$ can be replaced by $g+y_ic$, where $i=1,2$. For this, we first replace $g$ by $g+y_3c$ using $(b)$, i.e., we can assume that $g=G(y_1,y_2)+y_3c$. If $\g=(x_1,x_2,x_3+x_1,x_4,\ldots,x_n)$, then 
\bes
A(h,g)^{\g}=\s A(h^{\g},g^{\g}), 
\ees
where 
\bes
\s=(x_1,x_2,x_3-[x_1,x_n]h^{\g}g^{\g}-[x_2,x_n]h^{\g}(g^{\g})^2,x_4,\ldots,x_n)\\
=(x_1,x_2,x_3-[x_1,x_n]h^{\g}g^{\g},x_4,\ldots,x_n)(x_1,x_2,x_3-[x_2,x_n]h^{\g}(g^{\g})^2,x_4,\ldots,x_n)
\ees
is an almost tame automorphism. Consequently, $A(h,g)$ is almost tame if and only if $A(h^{\g},g^{\g})$ is almost tame. 
Note that 
\bes
g^{\g}=G(y_1,y_2)+y_3c+y_1c. 
\ees
Using $(b)$, again we can replace $g^{\g}$ by 
\bes
 G(y_1,y_2)+y_1c. 
\ees
This means that we can replace $g=G(y_1,y_2)$ by $G(y_1,y_2)+y_1c$ for any $c\in K[y_1,y_2]$. 

If $\delta=(x_1,x_2,x_3+x_2,x_4,\ldots,x_n)$ and $g=G(y_1,y_2)+y_3c$, then 
\bes
A(h,g)^{\delta}=\t A(h^{\delta},g^{\delta}), 
\ees
where 
\bes
\t=(x_1,x_2,x_3+[x_1,x_n]h^{\delta}+[x_2,x_n]h^{\delta}g^{\delta},x_4,\ldots,x_n)\\
=(x_1,x_2,x_3+[x_1,x_n]h^{\delta},x_4,\ldots,x_n)(x_1,x_2,x_3+[x_2,x_n]h^{\delta}g^{\delta},x_4,\ldots,x_n)
\ees
is again almost tame automorphism. In this case we have  
\bes
g^{\delta}=G(y_1,y_2)+y_3c+y_2c. 
\ees
Using $(b)$, we can replace $g^{\delta}$ by 
\bes
 G(y_1,y_2)+y_2c. 
\ees
This means that we can replace $g=G(y_1,y_2)$ by $G(y_1,y_2)+y_2c$ for any $c\in K[y_1,y_2]$. 

Thus  we can make $g=0$. Observe that 
\bes
A(h,0)=(x_1, x_2-[x_1,x_n]h,x_3,\ldots,x_n)
\ees
is a Chein automorphism.  $\Box$

\begin{theorem}\label{t3} Every automorphism of the form 
\bes
B(h,f,g)=(x_1+[x_1,x_n]hfg+[x_2,x_n]hg^2, x_2-[x_1,x_n]hf^2-[x_2,x_n]hfg,x_3,\ldots,x_n)
\ees
is almost tame.
\end{theorem}
\Proof  Consider the Chein automorphisms 
\bes
\phi_1=(x_1,x_2-[x_1,x_n]hf^2-[x_3,x_n]hf^2g,x_3,\ldots,x_n),\\
\psi_1=(x_1+[x_2,x_n]hg^2-[x_3,x_n]hfg^2,x_2,\ldots,x_n)
\ees 
and the automorphisms
\bes
\phi_2=(x_1+[x_1,x_n]hfg+[x_3,x_n]hfg^2, x_2,x_3-[x_1,x_n]hf-[x_3,x_n]hfg,x_4,\ldots,x_n),\\
\psi_2=(x_1,x_2-[x_2,x_n]hfg+[x_3,x_n]hf^2g,x_3-[x_2,x_n]hg+[x_3,x_n]hfg,x_4,\ldots,x_n). 
\ees
Note that 
\bes
(\phi_2)^{(23)}=A((hf)^{(23)},g^{(23)}) \ \hbox{and}  \ (\psi)^{132)}=A((hg)^{(132)},-f^{(132)}) 
\ees
are almost tame by Lemma \ref{l7}. Then $\phi_2$ and $\psi_2$ are also almost tame. Conseqently, 
\bes
\phi=\phi_1\phi_2=(x_1+[x_1,x_n]hfg+[x_3,x_n]hfg^2,x_2-[x_1,x_n]hf^2-[x_3,x_n]hf^2g,\\
x_3-[x_1,x_n]hf-[x_3,x_n]hfg,x_4,\ldots,x_n)
\ees
and 
\bes
\psi=(x_1+[x_2,x_n]hg^2-[x_3,x_n]hfg^2,x_2-[x_2,x_n]hfg+[x_3,x_n]hf^2g,\\
x_3-[x_2,x_n]hg+[x_3,x_n]hfg,x_4,\ldots,x_n)
\ees
are almost tame. Direct calculations give that  
\bes
\phi\psi=(x_1+[x_1,x_n]hfg+[x_2,x_n]hg^2, x_2-[x_1,x_n]hf^2-[x_2,x_n]hfg,\\
x_3-[x_1,x_n]hf-[x_2,x_n]hg,x_4,\ldots,x_n).
\ees
Then 
\bes
B(h,f,g)=\tau\phi\psi, 
\ees
where 
\bes
\tau=(x_1,x_2,x_3+[x_1,x_n]hf+[x_2,x_n]hg,x_4,\ldots,x_n).
\ees
Note that $\tau$ is a Chein automorphism. Hence $B(h,f,g)$ is almost tame. 
$\Box$

\section*{Acknowledgments}

The author would like to thank the Max Planck Institute f\"ur Mathematik  for
its hospitality and excellent working conditions, where some part of this work has been done. The author thanks Professors V. Dotsenko, M. Jibladze, L. Makar-Limanov, and I. Shestakov for useful discussions.

The project is supported by the grant AP14872073 of the Ministry of Education and Science of the Republic of Kazakhstan.


\begin{thebibliography}{99}




\bibitem{Artamonov} V.A. Artamonov, The categories of free metabelian groups and Lie algebras. 
Comment. Math. Univ. Carolinae 18 (1977), no. 1, 143--159.

\bibitem{Bachmuth} S. Bachmuth, 
Automorphisms of free metabelian groups. 
Trans. Amer. Math. Soc. 118 (1965), 93--104.

\bibitem{BM67} S. Bachmuth, H. Mochizuki, Automorphisms of a class of metabelian groups. II 
Trans. Amer. Math. Soc. 127 (1967), 294--301.

\bibitem{BM85} S. Bachmuth, H. Mochizuki, $\mathrm{Aut} F\to \mathrm{Aut}(F/F")$ is surjective for free group  $F$  of rank  $\geq 4$. 
Trans. Amer. Math. Soc. 292 (1985), no. 1, 81--101.


\bibitem{BN} Y. Bahturin, S. Nabiyev, Automorphisms and derivations of abelian extensions of some Lie algebras.
Abh. Math. Sem. Univ. Hamburg 62 (1992), 43--57.



\bibitem{Bahturin} Yu. Bahturin, Identical relations in Lie algebras. 
De Gruyter Exp. Math., 68
De Gruyter, Berlin, 2021, xxv+514 pp.



\bibitem{BD93-2} R.M. Bryant, V. Drensky, Dense subgroups of the automorphism groups of free algebras. 
Canad. J. Math. 45 (1993), no. 6, 1135--1154.

\bibitem{Chein} O. Chein, $\mathrm{IA}$  automorphisms of free and free metabelian groups. 
Comm. Pure Appl. Math. 21 (1968), 605--629.

\bibitem{Chibrikov} E. Chibrikov, 
On free Sabinin algebras. Comm. Algebra 39 (2011), no. 11, 4014--4035.


\bibitem{Cohn64}
P.M. Cohn, Subalgebras of free associative algebras. Proc. London Math. Soc. (3)14 (1964), 618--632.




\bibitem{Cohn06} P.M. Cohn, Free ideal rings and localization in general rings. 
New Math. Monogr., 3
Cambridge University Press, Cambridge, 2006, xxii+572 pp.




\bibitem{Czer}
A.J. Czerniakiewicz,  Automorphisms of a free associative algebra of
rank 2. I, II. Trans. Amer. Math. Soc. 160 (1971), 393--401; 171 (1972), 309--315.


\bibitem{DU22} V. Dotsenko, U. Umirbaev, An effective criterion for Nielsen-Schreier varieties, 
Int. Math. Res. Not. IMRN (2023), no. 23, 20385--20432.


\bibitem{Drensky92} V. Drensky, Wild automorphisms of nilpotent-by-abelian Lie algebras. 
Manuscripta Math. 74 (1992), no. 2, 133--141.

\bibitem{GGR} C.K. Gupta, N.D. Gupta, V.A. Roman'kov, Primitivity in free groups and free metabelian groups. 
Canad. J. Math. 44 (1992), no. 3, 516--523.


\bibitem{Jung}
H.W.E. Jung, Uber ganze birationale Transformationen der Ebene. J. reine angew Math. 184 (1942), 161--174. 





\bibitem{KR} A.N. Kabanov, V.A. Roman'kov, Strictly nontame primitive elements of a free metablelian Lie algebra of rank $3$. 
Sibirsk. Mat. Zh. 50 (2009), no. 1, 82--95.
Sib. Math. J. 50 (2009), no. 1, 66--76.



\bibitem{KMLU}
 D. Kozybaev, L. Makar-Limanov, U. Umirbaev,  The Freiheitssatz and the automorphisms of free right-symmetric algebras. 
Asian-Eur. J. Math. 1 (2008), no. 2, 243--254.

\bibitem{KSh} A.N. Krasil'nikov, A.L. Shmel'kin, 
Applications of the Magnus embedding in the theory of varieties of groups and Lie algebras. Fundam. Prikl. Mat. 5 (1999),  no. 2, 493--502.

\bibitem{Kulk}
W. van der Kulk, On polynomial rings in two variables. Nieuw Archief voor Wisk. (3)1 (1953), 33--41.

\bibitem{Kurosh} A.G. Kurosh, Non-associative free algebras and free products of algebras. Rec.Math. [Mat. Sbornik] N.S. 20/62 (1947), 239--262.

\bibitem{Lewin} J. Lewin,  On Schreier varities of linear algebras. Trans. Amer. Math. Soc. 132 (1968), 553--562.

\bibitem{Magnus} W. Magnus, On a theorem of Marshall Hall. 
Ann. of Math. (2) 40 (1939), 764--768.

\bibitem{ML70}
L.G. Makar-Limanov, The automorphisms of the free algebra of two
generators. Funktsional. Anal. i Prilozhen. 4  (1970), no. 3, 107--108; English translation: in Functional Anal. Appl. 4 (1970), 262--263.


\bibitem{MLSh}
L. Makar-Limanov, I. Shestakov, 
Polynomial and Poisson dependence in free Poisson algebras and free Poisson fields.
Journal of Algebra 349 (2012), 372--379.



\bibitem{MLTU}
L. Makar-Limanov, U. Turusbekova, U. Umirbaev,   Automorphisms 
and derivations of free Poisson algebras in two variables. J. Algebra 322 (2009), 3318--3330.

\bibitem{MLU16} L. Makar-Limanov, U. Umirbaev, Free Poisson fields
and their automorphisms. J. Algebra Appl. 15 (2016), no. 10, 1650196, 13 pp.



\bibitem{Mikhalev85} A.A. Mikhalev, Subalgebras of free colored Lie superalgebras. Mat. Zametki 37 (1985), no. 5, 653--661.

\bibitem{Mikhalev88} A.A. Mikhalev, Subalgebras of free Lie $p$-superalgebras. (Russian) Mat. Zametki 43 (1988), no. 2, 178--191, 300; translation in Math. Notes 43 (1988), no. 1--2, 99--106.

\bibitem{MS14} A.A. Mikhalev, I.P. Shestakov, 
PBW-pairs of varieties of linear algebras. Comm. Algebra 42(2014), no.2, 667--687.


\bibitem{Nagata} M. Nagata, On the automorphism group of $k[x,y]$. Lect. in Math., Kyoto Univ., Kinokuniya, Tokio, 1972.



\bibitem{Nauryzbaev09} R. Nauryzbaev, Reducibility of automorphisms of free metabelian Lie algebras of rank $3$. (Russian) Herald of the Eurasian National University, (2009), no. 6, 200--213. 

\bibitem{Nauryzbaev10} R. Nauryzbaev, Defining relations of the tame automorphism group of free metabelian Lie algebras of rank $3$. (Russian) Herald of the Eurasian National University, (2010), no. 4, 164--170. 

\bibitem{Nauryzbaev11} R. Nauryzbaev, Defining relations of the almost tame automorphism group of free metabelian Lie algebras of rank $3$. (Russian) Herald of the Eurasian National University, (2011), no. 4, 21--26. 



\bibitem{Nauryzbaev15} R. Nauryzbaev, On generators of the tame automorphism group of free metabelian Lie algebras. 
Comm. Algebra 43 (2015), no. 5, 1791--1801.

\bibitem{NSU25} R. Nauryzbaev, I. Shestakov, U. Umirbaev, Chein automorphisms of free metabelian anticommutative algebras. Serdica Math. J. 51 (2025), 351--366. 


\bibitem{Nielsen21} J. Nielsen, Om regning med ikke-kommutative faktorer og dens anvendelse i gruppeteorien, Math. Tidsskrift B (in Danish), 1921, 78--94. 

\bibitem{Nielsen24} J. Nielsen, Die Isomorphismengruppe der freien Gruppen, Mathematische Annalen (in German), 91 (3--4) (1924), 169--209.

\bibitem{OE} Z. \"Ozkurt, N. Ekici, An application of the Dieudonné determinant: detecting non-tame automorphisms. 
J. Lie Theory 18 (2008), no. 1, 205--214.

\bibitem{Romankov85} V.A. Roman'kov, The automorphism group of a free metabelian group. The Interconnection Problems of Abstract and Applied Algebra, Novosibirsk (1985), pp. 53--80 (in Russian).

\bibitem{Romankov92} V.A. Roman'kov, Primitive elements of free groups of rank  $3$. 
Math. USSR-Sb. 73 (1992), no. 2, 445--454.

\bibitem{Romankov95} V.A. Roman'kov, The Tennant-Turner swap conjecture. 
Algebra i Logika 34 (1995), no. 4, 448--463, 487--488.
Algebra and Logic 34 (1995), no. 4, 249--257.

\bibitem{Romankov08} V.A. Roman'kov, On the automorphism group of a free metabelian Lie algebra
Internat. J. Algebra Comput. 18 (2008), no. 2, 209--226.

\bibitem{Romankov20} V.A. Roman'kov, Primitive elements and automorphisms of the free metabelian group of rank 3. 
Sib. Elektron. Mat. Izv. 17 (2020), 61--76.





\bibitem{SU02}
I.P. Shestakov and U.U. Umirbaev, Free Akivis algebras, primitive elements and their
relations with loops. J. Algebra 250 (2002), 533--548.

\bibitem{SU04-1}
I.P. Shestakov, U.U. Umirbaev, Tame and wild automorphisms of rings of polynomials in three variables. J. Amer. Math. Soc. 17 (2004), 197--227.

\bibitem{SU25} Shestakov, I., Umirbaev, U.: Tangent Lie algebras of automorphism groups of free algebras. Linear Algebra Appl. (accepted) ArXiv. https://arxiv.org/abs/2507.20486 (2025). 

\bibitem{SZ24} I. Shestakov, Z. Zhang, An Anick type wild automorphism of free Poisson algebras. ArXiv. https://arxiv.org/abs/2407.04919 (2024), 15 pages. 

\bibitem{Shirshov53} A.I. Shirshov, Subalgebras of free Lie algebras. (Russian) Mat. Sbornik 33/75 (1953), 441--452.


\bibitem{Shirshov54} A.I. Shirshov,  Subalgebras of free commutative and free anticommutative algebras. (Russian) Mat. Sbornik 34/76 (1954), 81--88.




\bibitem{Shmelkin73} A.L. Shmel'kin, Wreath products of Lie algebras and their application in the theory of groups, Trudy 
Moskov. Mat. Obshch. 29 (1973), 247--260, Russian; Translation, Trans. Moscow Math. Soc. 29 (1973), 
239--252. 

\bibitem{Schreier} O. Schreier, 
Die Untergruppen der freien Gruppen. (German) Abh. Math. Sem. Univ. Hamburg 5(1927), no.1, 161--183.

\bibitem{Shpilrain} V. Shpilrain, On generators of  $L/R^2$  Lie algebras. 
Proc. Amer. Math. Soc. 119 (1993), no. 4, 1039--1043.

\bibitem{Stern}  A.S. Shtern,  Free Lie superalgebras. Sibirsk. Mat. Zh. 27 (1986),  no. 1, 170--174.

\bibitem{Smith} M.K. Smith, Stably tame automorphisms,
J. Pure and Appl. Algebra, 58(1989), 209--212.



\bibitem{U93} U.U. Umirbaev,  Partial derivations and endomorphisms of some relatively free Lie algebras. (Russian) Sibirsk. Mat. Zh. 34 (1993), no. 6, 179--188. 

\bibitem{U94}
U.U. Umirbaev, On Schreier varieties of algebras. Algebra i logika  33 (1994), no. 3, 
317--340; Ehglish translation: In Algebra and Logic  33 (1994), no. 3, 180--193. 

\bibitem{95SMJ}
 U.U. Umirbaev, On an extension of automorphisms of polynomial rings. (Russian) Sibirsk. Mat. Zh. 36 (1995), no. 4, 911--916; translation in Siberian Math. J. 36 (1995), no. 4, 787--791.

\bibitem{U96} U.U. Umirbaev, Universal derivations and subalgebras of free algebras.
 Proc. 3rd Internat.
    Conf. Algebra, Krasnoyarsk, Russia,
    Walter de Gruyter, Berlin, 1996, 255--271.

\bibitem{US02RAN} U.U. Umirbaev, I.P. Shestakov,  Subalgebras and automorphisms of polynomial rings. (Russian) Dokl. Akad. Nauk 386 (2002), no. 6, 745--748.

\bibitem{UU06AN} U.U. Umirbaev, Defining relations of the tame automorphism group of polynomial rings, and wild automorphisms of free associative algebras. (Russian) Dokl. Akad. Nauk 407 (2006), no. 3, 319--324.


\bibitem{UU06JR}
 U.U. Umirbaev,  Defining relations of the tame automorphism group of
 polynomial algebras in three variables.  J. Reine Angew. Math. 600 (2006), 203--235.

\bibitem{UU07JR}
 U. Umirbaev,  The Anick automorphism of free associative algebras.  J. Reine Angew. Math. 605 (2007), 165--178.


\bibitem{UU07JA}
U.U. Umirbaev, Defining relations for automorphism groups of free algebras. J. Algebra 314 (2007), no. 1, 209--225.





\bibitem{U12}
U. Umirbaev, Universal enveloping algebras and universal derivations
of Poisson algebras. J. Algebra 354 (2012), pp. 77--94.

\bibitem{U24} U. Umirbaev, Erroneous proofs of the wildness of some automorphisms of free metabelian Lie algebras. ArXiv. /abs/2401.07182, (2024), 10 pages, https://doi.org/10.48550/arXiv.2401.07182



\bibitem{Witt} E. Witt, Die Unterringe der freien Lieschen Ringe. Math. Z.  64 (1956), 195--216.


\end{thebibliography}
\end{document}